\documentclass[a4paper, 11pt]{article}

\usepackage{amsmath,amsfonts,amssymb,amsthm}
\usepackage[a4paper, margin=1in]{geometry}
\usepackage{enumitem}

\usepackage[T1]{fontenc}

\usepackage{lmodern}
\usepackage{calc}
\usepackage[dvipsnames]{xcolor}
\usepackage{nameref}
\usepackage{upgreek}     
\usepackage{enumitem}
\usepackage{graphicx}
\usepackage{mathrsfs}
\usepackage{tikz-cd}
\usepackage{hhline}
\usepackage{textcomp}
\usepackage{amscd}
\usepackage[all,cmtip]{xy}
\usepackage{mathtools}
\usepackage{extarrows}
\usepackage{upgreek}
\usepackage{bbm}
\usepackage[displaymath]{lineno}
\usepackage[singlespacing]{setspace}%
\usepackage[colorlinks=true,breaklinks=true,bookmarks=true,urlcolor=blue,
citecolor=BrickRed,linkcolor=teal,bookmarksopen=false,draft=false]{hyperref}

\providecommand{\U}[1]{\protect\rule{.1in}{.1in}}

\makeatletter
\let\orgdescriptionlabel\descriptionlabel
\renewcommand*{\descriptionlabel}[1]{%
	\let\orglabel\label
	\let\label\@gobble
	\phantomsection
	\edef\@currentlabel{#1}%
	\let\label\orglabel
	\orgdescriptionlabel{#1}%
}
\makeatother

\theoremstyle{plain}
\newtheorem{theorem}{Theorem}[section]

\newtheorem{corollary}[theorem]{Corollary}

\newtheorem{remark}[theorem]{Remark}
\theoremstyle{definition}
\newtheorem{definition}[theorem]{Definition}
\newtheorem*{definition*}{Definition} 

\newcommand{\RR}{\mathbb{R}}

\newcommand{\ps}{\mathrm{({PS}})_m}

\sbox0{$\fontdimen8\textfont3=0.4pt$}

\renewcommand{\emptyset}{\varnothing}

\newcommand{\abs}[1]{\lvert#1\rvert}

\DeclareMathAlphabet{\mathpzc}{OT1}{pzc}{m}{it}

\newcommand{\ftt}[1] {\mathsf{#1}}

\newcommand{\di}{differentiable }
\newcommand{\ga}{G\^{a}teaux }
\newcommand{\va}{\varphi}

\newcommand{\cs}{continuous }

\newcommand{\dd}{{\tt D}}

\newcommand{\fs}[1]{\mathsf {#1}}
\newcommand{\eu}[1]{\EuScript {#1}}

\newcommand*{\medcap}{\mathbin{\scalebox{1.5}{\ensuremath{\cap}}}}

\newcommand{\mt}{\mathbbm {d}}

\newcommand{\set}[1]{\left\{#1\right\}}
\newcommand\Set[2]{\left\{#1\mid#2\right\}} 
\newcommand{\snorm}[2][]{\left\lVert#2\right\rVert_{#1}}

\newcommand{\Sem}[1]  {\textsf{Sem}(#1)}

\newcommand{\zero}[1]{\boldsymbol{0}_{#1}}

\usepackage{bm}    

\newcommand{\rr}{\mathbb{R}}
\newcommand{\nn}{\mathbb{N}}

\DeclareMathAlphabet\EuScript{U}{eus}{m}{n}
\SetMathAlphabet\EuScript{bold}{U}{eus}{b}{n}
\newcommand\opn{\ensuremath{\mathrel{\mathpalette\opncls\circ}}}

\newcommand{\opncls}[2]{
	\ooalign{$#1\subseteq$\cr
		\hidewidth\raisefix{#1}\hbox{$#1{\stylefix{#1}#2}\mkern2mu$}\cr}}
\def\raisefix#1{
	\ifx#1\displaystyle
	\raise.39ex
	\else
	\ifx#1\textstyle
	\raise.39ex
	\else
	\ifx#1\scriptstyle
	\raise.275ex
	\else
	\raise.150ex
	\fi
	\fi
	\fi
}
\def\stylefix#1{
	\ifx#1\displaystyle
	\scriptstyle
	\else
	\ifx#1\textstyle
	\scriptstyle
	\else
	\ifx#1\scriptstyle
	\scriptscriptstyle
	\else
	\scriptscriptstyle
	\fi
	\fi
	\fi
}
\DeclareFontFamily{U}{mathx}{\hyphenchar\font45}
\DeclareFontShape{U}{mathx}{m}{n}{
	<5> <6> <7> <8> <9> <10>
	<10.95> <12> <14.4> <17.28> <20.74> <24.88>
	mathx10
}{}
\newcommand{\fr}{Fr\'{e}chet }

\newcommand{\bo}{\mathbf{B}^{\mathrm{op}}}

\DeclareMathAlphabet{\mathsfit}{OT1}{cmss}{m}{sl}

\newcommand{\Cl}[1]{\overline{#1}}
\newcommand{\subjclass}[1]{\textbf{AMS Subject Classifications (2010):} #1\par}
\newcommand{\keywords}[1]{\textbf{Keywords:} #1\par}
\title{On the existence of a global diffeomorphism between Fr\'{e}chet spaces}
\author{Kaveh Eftekharinasab}
\date{}
\usepackage{cite}
\begin{document}
	
	\maketitle
	
	\begin{abstract}
	We provide sufficient conditions for the existence of a global diffeomorphism between tame Fr\'{e}chet spaces.
We prove a version of the Mountain Pass Theorem which is a key ingredient
in the proof of the main theorem. 
	\end{abstract}
\bigskip

	\let\thefootnote\relax\footnotetext{
		\,\,\,\,\,\keywords{The Nash-Moser inverse function theorem, Mountain Pass Theorem, Ekelend's variational principle, Keller's \(C_c^k\)-mappings.}}
	\section*{Introduction}
In this paper we consider the problem of finding sufficient conditions under which a tame map
between tame Fr\'{e}chet spaces becomes a global diffeomorphism. Tame maps are important because they appear not only as
differential equations but also as their solutions (see~\cite{ham} for examples).
Although the theory of differential equations in Fréchet spaces has a significant relation with problems in both linear and nonlinear functional analysis, not many methods for solving different types of differential equations are known. 
Our result would provide an approach to solve an initial value nonlinear  integro-differential equation 
$$
x'(t) + \int_0^t \phi (t,s,x(s)) ds = y(t) \quad t \in [0,1].
$$
We follow the ideas in~\cite{I} and~\cite{G} where the analogue problem for Banach and Hilbert spaces was studied.  
There are two approaches to calculus on Fréchet spaces: the Gâteaux-approach (see~\cite{ke}) and the so-called convenient analysis (see~\cite{kri}). We will apply the first one because, to define the Palais-Smale condition, which plays a significant role in the calculus of variation, we need an appropriate topology on dual spaces that is compatible with our notion of differentiability; only in the first approach does such a topology exist.

In~\cite{k2}, the author defined the Palais-Smale condition for Keller's $C_c^1$-maps between Fr\'{e}chet spaces and obtained some existence results for locating critical points. In this paper, by means of this condition we generalize the  mountain pass theorem of Ambrosetti and Rabinowitz to Fr\'{e}chet spaces. Our proof of the mountain pass theorem relies on the Ekeland's variational principle. Since, in general, we cannot acquire deformation results for Fr\'{e}chet spaces due to the lack of a general solvability theory for differential equations, our proof of the mountain pass theorem relies on the Ekeland's variational principle. It is worth mentioning that for every Fr\'{e}chet space the projective limit technique gives a way to solve a wide class of differential equations (see~\cite{dod}). This technique would be a way of obtaining many results such as deformation lemmas.

Roughly speaking, the main theorem states that if $ \varphi $
is a smooth tame map that satisfies the assumptions of the Nash-Moser inverse function theorem, and if, for an appropriate  auxiliary functional
$ \iota $, the functional
$e \mapsto \iota (\varphi (e) - f) $  satisfies the Palais-Smale condition at any level, then $ \varphi $ is a global diffeomorphism. 
\section{Mountain Pass Theorem}

Throughout this paper, we assume that \( (\fs{F}, \Sem{\fs{F}}) \) and \( (\fs{E}, \Sem{\fs{E}}) \) are \fr spaces over \( \rr \), where \( \Sem{\fs{F}} = \Set{\snorm[\fs{F},n]{\cdot}}{n \in \nn} \) and \( \Sem{\fs{E}} = \Set{\snorm[\fs{E},n]{\cdot}}{n \in \nn} \) are increasing sequences of continuous seminorms that define the topologies of \( \fs{F} \) and \( \fs{E} \), respectively. 
We also apply the translation-invariant metric $\mt$ given by
\[
\mt(x,y) \coloneqq \sum_{n=1}^{\infty} \dfrac{1}{2^n}\dfrac{\snorm[\fs{F},n]{x-y}}{1+\snorm[\fs{F},n]{x-y}},
\]
which generates a topology equivalent to the one defined by the family of seminorms.

We use the notation \( U \opn \mathsf{T} \) to denote that \( U \) is an open subset of the topological space \( \mathsf{T} \).

Let \( \mathfrak{S} \) be a family of bounded subsets of \( \fs{F} \), with the following two properties:
\begin{itemize}
	\item[(\(\mathfrak{S}_1\))] If \( A, B \in \mathfrak{S} \), then there exists \( C \in \mathfrak{S} \) such that \( A \cup B \subset C \).
	\item[(\(\mathfrak{S}_2\))] If \( A \in \mathfrak{S} \) and \( r \) is a real number, then there exists \( B \in \mathfrak{S} \) such that \( rA \subset B \).
\end{itemize}

Throughout the paper, we assume that \( \mathfrak{S} \) is the family of all compact subsets of \( \fs{F} \). We consider the topology of compact convergence on the dual space, and denote the dual space endowed with this topology with by \( \fs{F}'_c \).

 Let \( \mathcal{L}_c(\fs{F}, \fs{E}) \) be the space of all continuous linear mappings from \( \fs{F} \) to \( \fs{E} \) endowed with the topology of compact convergence, which is a Hausdorff locally convex topology, defined by seminorms:
 \[
 \snorm[{S},i]{\ell} \coloneqq \sup\Set{\snorm[\fs{E},i]{\ell (f)}}{ f  \in S}
\]
where \(S \in \mathfrak{S}\). If \(\fs{E}=\RR \) with the usual modulus \(\abs{\cdot}\), then
we denote by \(\snorm[{S}]{\cdot}\) the seminorms that define the topology of \( \mathcal{L}_c(\fs{F}, \rr) \) 
\begin{definition}[Definition 1.0.0, \cite{ke}]\label{def:diff}
	Let  $ \va\colon U \opn \fs{E}  \to  \fs{F}$ be a mapping. Then the  derivative
	of $\va$ at $x$ in the direction $h$ is defined by 
	\[
	\dd \va_x(h)=\dd\va(x)(h) \coloneqq
	\lim_{t \to 0} {1\over t}(\va(x+th) -\va(x))
	\]
	whenever it exists. 
	The function $\va$ is called differentiable at
	$x$ if $\dd \va(x)(h)$ exists for all $h \in \fs{E}$. 
	It is called \(C_c^1\)-mapping
	if it is differentiable at all
	points of $U$, and the mapping
	\[
	\dd \va \colon U  \to  \mathcal{L}_c(\fs{E}, \fs{F}) 
	\]
	is continuous. 
\end{definition}
 Higher order differentiability is defined in \cite[Definition 2.5.0]{ke}.
The primary motivation for employing this class of mappings is the need to equip dual spaces with a suitable topology in order to define the Palais–Smale condition. These mappings are known as Keller \( C_c^k \)-mappings, which are equivalent to  widely used Michal--Bastiani notion of differentiability.

\begin{definition}[The Palais-Smale Condition, Definition 1.2, \cite{k2}]\label{def:PS}
	Let $\va\colon U \opn \fs{F} \to \rr$ be a Keller $C_c^1$-functional.
	\begin{itemize}
		\item [(i)] We say that $\va$ satisfies the Palais-Smale condition, $\mathrm{({PS}})$-condition in short, if every sequence $(x_i) \subset \fs{F}$ such that $\va(x_i)$ is bounded and
		$$
		\dd \va(x_i) \to 0 \quad \text{in} \quad \fs{F}_{c}',
		$$
		has a convergent subsequence.
		\item[(ii)] We say that $\va$ satisfies the Palais-Smale condition at level $m \in \rr$, $\mathrm{({PS}})_m$-condition in short, if every sequence $(x_i) \subset \fs{F}$ such that
		$$
		\va(x_i) \to m \quad \text{and} \quad \dd \va(x_i) \to 0 \quad \text{in} \quad \fs{F}_{c}',
		$$
		has a convergent subsequence.
	\end{itemize}
\end{definition}

If $ \va $ satisfies the $ \ps $-condition, then every $ \ps $-sequence converges, up to a subsequence,
to some point $ p $  and, by continuity, one has that  $ \va(p)=m $ and $ \va'(p) =0 $. In another word, $ p $
is a critical point of $ \va $.

\begin{corollary}[Corollary 4.9, \cite{k2}]\label{co:minimizing}
	Let  $\va\colon\fs{F}\to \rr$ be a Keller's \(C_c^1\)-mapping which is bounded from below. 
	If the $(\mathrm{PS})_m$-condition holds with $m = \inf_{\fs{F}} \va$, then $\va$ attains its minimum at a critical point $ x_0 \in \fs{F} $ with $ \va(x_0) =m $.
\end{corollary}

Consider the following weak form of the Ekeland's variational principle.

\begin{theorem}[Theorem 1 bis, \citen{mini}]\label{th:weke}
	Assume that $(\mathsf{M}, \mathbbm{m})$ is a complete metric space.
	Let a functional $\upphi \colon \mathsf{M} \to (-\infty, \infty]$ be lower semicontinuous, bounded from below, and not identically equal to $\infty$.
		Then, for any $ \epsilon > 0 $ there exists $m \in \fs{M}$ such that
	\begin{enumerate}
		\item $\upphi(m) \le \inf_{ \fs{M}} + \epsilon $;
		\item $\upphi(m) < \upphi(m) + {\epsilon} \mathbbm{m} (m,n) \quad \forall n \in \mathsf{M} \setminus \{m\}$.
	\end{enumerate}
\end{theorem}

The mountain pass theorem  is a minimax result.
This theorem relies on a geometric condition on functionals, where a relation exists between the values of the functional over minimax sets. This condition is described as follows: let $ U $ be an open set such that $ x_0 \in U$ and $x_1 \notin \Cl{U}$, such that
\[
\max \set{\va(x_0), \va (x_1)} \leq \inf_{\partial U} \va.
\]

Next, we will describe the choice of class of sets in the minimax expression. 
We denote by $ \ftt{C}^1(\fs{F},\rr) $ the set of Keller $ C_c^1 $-functionals. Let $ \va \in \ftt{C}^1(\fs{F},\rr) $ be given, and
let 
\[
\Gamma \coloneqq \Big\{ \gamma \in \ftt{C} ([0,1], \fs{F}) : \gamma (0)= x_0,\gamma(1)=x_1\in \fs{F} \Big\}
\]
 be the set of continuous paths joining $ x_0 $ and $ x_1 $. Consider the Fr\'{e}chet space $\ftt{C} ([0,1], \fs{F})$ with the family of seminorms: 
\begin{equation*}
	\varphi \mapsto \snorm[\ftt{C},n]{\va} \coloneqq \sup_{0 \leq t \leq 1} \snorm[\fs{F},n]{\va(t)}.
\end{equation*}
The metric 
\begin{equation}\label{met:dfg}
	\mt_{\ftt{C}} (\gamma,\eta) \coloneqq \sum _{n \in \nn} \dfrac{\snorm[\ftt{C},n]{\gamma - \eta}}{1+\snorm[\ftt{C},n]{\gamma - \eta}}
\end{equation}
is complete, translation-invariant, and induces the same topology on  $ \ftt{C} ([0,1], \fs{F}) $.
We can easily show that $\Gamma$ is closed in $\ftt{C} ([0,1], \fs{F})$ and thus it is a complete metric space with the metric $\mt_{\Gamma}$, which is the restriction of $\mt_{\ftt{C}}$ to $\Gamma$. Then, we consider the minimax expression of the following form:
\[
c \coloneqq \inf_{\gamma \in \Gamma} \max_{t \in [0,1]} \va(\gamma(t)). 
\]

The idea of the proof of the mountain pass theorem is straightforward: for a given  $ \va \in \ftt{C}^1(\fs{F},\rr) $ that satisfies the PS-condition and a point $x_0 \in \fs{F}$, if a particular condition holds (the condition~\eqref{mpt:geo}), we define a functional $\Psi$ on $\Gamma$ so that it satisfies the assumptions of the Ekeland variational principle (Theorem~\ref{th:weke}).
Then this theorem yields that $\Psi$ has almost minimizer points satisfying certain conditions. We use a
sequence of these points on $\Gamma$ and associate this sequence of almost minimizers  with a sequence on $\fs{F}$, which satisfies the requirements of 
the  PS-condition for $ \va $. The limit of a subsequence of this sequence in $\fs{F}$ is a critical point
of $ \va $. The difficult step is to find a connection between the sequence of almost minimizers of $\Psi$ and a sequence in $\fs{F}$ that satisfies the PS-condition.
\begin{theorem}[ Theorem 1, \cite{mean}]\label{th:lmvt}
	Let $ \va : A \opn \fs{F} \to \rr $ be a \cs and \ga \di function
	at each point of $ [x_0,x_0+h] $ in $ A $. Then there  exists a $ r \in (0,1) $ such that
	\begin{equation*}
		\va(x_0+h) - \va(x_0) = \va'(x_0+rh)(h).
	\end{equation*}
\end{theorem}
\begin{theorem}\label{th:mpt} Assume that $ \va \in \ftt{C}^1(\fs{F},\rr) $ satisfies the Palais-Smale condition at all levels. Let $ x_0 \in \fs{F} $ and suppose that $\va$ satisfies the condition\textup{:} 
	\begin{equation}\label{mpt:geo}
		\inf_{p \in \partial U} \va(p) > \max \{ \va(x_0),\va(x_1)\} = a,
	\end{equation}
	where $ \partial U$ is the boundary of an open neighborhood $U$ of $x_0 \in \fs{F}$ such that $x_1$  belongs to the distinct \textup{(}arcwise\textup{)} connected component  of $\fs{F} \setminus \partial{U} $.
	Then $ \va $ has a critical value $ c > a $ which can be characterized as
	\[
	c \coloneqq \inf_{\gamma \in \Gamma} \max_{t \in [0,1]} \va(\gamma(t)). 
	\]
\end{theorem}
\begin{proof}
The metric $\mt_{\Gamma}$, which is the restriction of the metric $\mt_{\ftt{C}}$ from~\eqref{met:dfg} to $\Gamma$, defines the topology of $\Gamma$. As established previously, $\Gamma$ is a complete metric space with this metric.

Define the functional $ \Psi : \Gamma \rightarrow \rr $ by
\[
\Psi(\gamma) = \max_{t \in [0,1]} \va(\gamma(t)).
\]
Since $\va$ is continuous (as a $C_c^1$-functional) and $\Psi$ is the supremum of a family of continuous functions over a compact set, it follows that $\Psi$ is lower semicontinuous.

Given that $x_0 \in U$ and $x_1 \notin \Cl{U}$, and $x_1$ lies in a distinct (arcwise) connected component of $\fs{F} \setminus \partial U$ from $x_0$, any continuous path $\gamma \in \Gamma$ from $x_0$ to $x_1$ must intersect the boundary $\partial U$. Thus, for all $\gamma \in \Gamma$, we have
\begin{equation}\label{mpt:6}
	\gamma ([0,1]) \medcap \partial U \neq \emptyset.
\end{equation}
Therefore, for any $\gamma \in \Gamma$, there exists $t_0 \in [0,1]$ such that $\gamma(t_0) \in \partial U$. This implies
\[
\max_{t \in [0,1]} \va(\gamma(t)) \geq \va(\gamma(t_0)) \geq \inf_{p \in \partial U} \va(p).
\]
Taking the infimum over all $\gamma \in \Gamma$, we get
\begin{equation}
	c = \inf_{\gamma \in \Gamma} \max_{t \in [0,1]} \va(\gamma(t)) \geq \inf_{p \in \partial U} \va(p) = c_1.
\end{equation}
From the geometric condition~\eqref{mpt:geo}, we have $ c_1 > a $.
Consequently, $c > a$. Since $\Psi(\gamma) \ge c$ for all $\gamma \in \Gamma$, and $c$ is a finite value, $\Psi$ is bounded from below.	

Let $ \widehat{\gamma}\in \Gamma $. We'll show that $ \Psi $ is continuous at $ \widehat{\gamma} $. Given $ \varepsilon > 0 $, by the continuity of $\va$ on the compact set $\widehat{\gamma}([0,1])$, we can choose $ \varrho > 0$ such that for all $y \in \widehat{\gamma}([0,1])$ and all $x \in \fs{F}$ with $\mt_{\fs{F}} (x,y) < \varrho $, we have $ \mid \va(x) - \va(y) \mid < \varepsilon$.
Then, for each $ \overline{\gamma} \in \Gamma$ such that $ \mt_{\Gamma}(\widehat{\gamma},\overline{\gamma}) < \varrho$, we have
\[
\Psi(\overline{\gamma}) - \Psi(\widehat{\gamma}) = \va (\overline{\gamma}(t_m))- \max_{t \in [0,1]}\va (\widehat{\gamma}(t)) \leq \va (\overline{\gamma}(t_m)) - \va (\widehat{\gamma}(t_m)),
\]
where $ t_m \in [0,1] $ is the point where the maximum of $ \va (\overline{\gamma}(t)) $ is attained. Since
\[
\mt_{\fs{F}} (\widehat{\gamma}(t_m), \overline{\gamma}(t_m) ) \leq \mt_{\Gamma}(\widehat{\gamma},\overline{\gamma}) < \varrho,
\]
it follows from our choice of $\varrho$ that $ \Psi(\overline{\gamma}) - \Psi(\widehat{\gamma}) < \varepsilon$. By symmetry (reverting the roles of $ \widehat{\gamma} $ and $ \overline{\gamma} $), we similarly find that $\Psi(\widehat{\gamma}) - \Psi(\overline{\gamma}) < \varepsilon$. Combining these, we obtain $ \mid \Psi(\widehat{\gamma})- \Psi(\overline{\gamma}) \mid < \varepsilon.$

Thus $ \Psi $ satisfies all conditions of Theorem~\ref{th:weke} (Ekeland's Variational Principle), and hence, for every $ \epsilon > 0 $ there exists $ \gamma_{\epsilon} \in \Gamma$ such that
\begin{align}
	\Psi (\gamma_{\epsilon}) 
	&\leq c + \epsilon, \label{mpt:5} \\
	\Psi (\gamma_{\epsilon}) 
	&\leq \Psi (\gamma) + \epsilon \, \mt_{\Gamma} (\gamma, \gamma_{\epsilon}),
	\quad \forall \gamma \neq \gamma_{\epsilon} \in \Gamma. \label{mpt:dddd}
\end{align}
Without loss of generality, we may assume
\begin{equation}\label{mpt:aa}
	0 < \epsilon < c - a.
\end{equation}
Our goal is to show that there is an $s \in [0,1]$ such that for every compact set $B \in \mathfrak{S}$ (defining the seminorms on $\fs{F}'_c$), we have
\begin{equation}\label{mpt:cont1}
	\snorm[B]{\va'(\gamma_{\epsilon}(s))}   \leq  \epsilon.
\end{equation}
We prove the inequality~\eqref{mpt:cont1} by contradiction. Recall that for any $B \in \mathfrak{S}$, the seminorm on $\fs{F}'_c$ is given by $\snorm[B]{\ell} = \sup_{f \in B} \abs{\langle \ell, f \rangle}$. Thus, for $\va'(\gamma_{\epsilon}(s))$, we have
\begin{equation*}
	\snorm[B]{\va'(\gamma_{\epsilon}(s))} = \sup_{g \in B} \abs{\langle \va'(\gamma_{\epsilon}(s)), g \rangle}.
\end{equation*}
Define the set
\begin{equation}
	S(\epsilon) \coloneqq \Set{s \in [0,1 ]}{c - \epsilon \leq \va(\gamma_{\epsilon}(s))}.
\end{equation}
From~\eqref{mpt:aa}, we have $a < c - \epsilon$. Since $\gamma_{\epsilon}(0)=x_0$ and $\va (x_0) \leq a$, it follows that $\va(\gamma_{\epsilon}(0)) \leq a < c - \epsilon$, which implies $0 \notin S(\epsilon)$. Furthermore, as $\va$ is continuous on $\fs{F}$ and $\gamma_{\epsilon}$ is continuous, the function $\va \circ \gamma_{\epsilon}$ is continuous. Therefore, $S(\epsilon)$ is the preimage of the closed set $[c-\epsilon, \infty)$ under a continuous map, making it a closed subset of the compact interval $[0,1]$. Thus, $S(\epsilon)$ is compact.

Suppose, to the contrary, that for all $s \in [0,1]$, the inequality~\eqref{mpt:cont1} does not hold. This means that for each $s \in S(\epsilon)$, there exists a compact set $B_s \in \mathfrak{S}$ such that
\[
\snorm[S_s]{\va'(\gamma_{\epsilon}(s))} > \epsilon
\]
 From the definition of the seminorm, this implies there exists an element $g_s \in B_s$ such that $\abs{\langle \va'(\gamma_{\epsilon}(s)), g_s \rangle} > \epsilon$. We can assume without loss of generality that $g_s$ is chosen such that 
\begin{equation}\label{mpt:impo}
	\langle \va'(\gamma_{\epsilon}(s)), g_s \rangle < - \epsilon.
\end{equation}
Since $\va'$ is continuous, it follows from~\eqref{mpt:impo} that for each $s \in S(\epsilon)$, there exist a scalar $\alpha_{s} > 0$ and an open interval $B_s \subset [0,1]$ containing $s$, such that for any $t \in B_s$ and any $h \in \fs{F}$ with $\mt_{\fs{F}}(\gamma_{\epsilon}(t)+h, \gamma_{\epsilon}(t)) < \alpha_s$ (i.e., $h$ is sufficiently small in the $\fs{F}$ metric), the inequality holds:
\begin{equation}\label{tmp:bb}
	\langle \va'(\gamma_{\epsilon}(t)+h), g_s \rangle < - \epsilon.
\end{equation}
The family $\{ B_s\}_{s \in S(\epsilon)}$ forms an open cover of the compact set $S(\epsilon)$, so there exists a finite subcovering $\{ B_{s_1}, \cdots, B_{s_k} \}$ of $S(\epsilon)$. Since $0 \notin S(\epsilon)$, we can assume $0 \notin B_{s_i}$ for any $i$. Consequently, for each $i=1,\dots,k$, $[0,1]\setminus B_{s_i}$ is a closed and non-empty set. 
Therefore, if $t \in \bigcup_{i=1}^k B_{s_i}$, then
\[
\sum_{i=1}^k\mathrm{dist} (t, [0,1]\setminus B_{s_i}) > 0.
\]
Now, we define functions
\(
\rho_j(t) = \mathrm{dist}(t, [0,1]\setminus B_{s_j})
\)
 for $j=1,\dots,k$. Note that $\rho_j(t) > 0$ if and only if $t \in B_{s_j}$. 

Then, define the functions $\chi_j(t) : [0,1] \rightarrow [0,1]$ by
\begin{gather*}
	\chi_j(t) =\begin{cases}
		\dfrac {\rho_j(t)}{\sum_{i=1}^k \rho_i(t)} & \text{if } t \in \bigcup_{i=1}^k B_{s_i}, \\
		0 & \text{otherwise}.
	\end{cases}
\end{gather*}
It is easily seen that each $\chi_j$ is continuous, $\chi_j(t) = 0$ if $t \notin B_{s_j}$, and importantly,
\begin{equation} \label{tmp:idc}
	\sum_{j=1}^k \chi_j(t) = 1 \quad \text{for } t \in \bigcup_{i=1}^k B_{s_i}.
\end{equation}
Fix a continuous function $\chi : [0,1] \to [0,1]$ such that
\begin{gather*}
	\chi(t) =\begin{cases}
		1 & \text{if } \va (\gamma_{\epsilon}(t)) \geq c, \\
		0 & \text{if } \va (\gamma_{\epsilon}(t)) \leq c - \epsilon.
	\end{cases}
\end{gather*}
The existence of such a continuous function is guaranteed by the continuity of $\va \circ \gamma_{\epsilon}$ and Urysohn's Lemma.
For this arbitrary fixed $\epsilon$, define the continuous function $ \mu : [0,1] \rightarrow \fs{F}$ by
$$
\mu(t)= \gamma_{\epsilon}(t)+ \alpha \chi (t) \sum_{j=1}^k \chi_j(t)g_{s_j}.
$$
Here $ \alpha = \min \{ \alpha_{s_1}, \cdots ,\alpha_{s_k}\} $, where $g_{s_j}$ are the vectors chosen in \eqref{mpt:impo}.
Now, we show that $\mu \in \Gamma$. By our choice of $\epsilon$ in~\eqref{mpt:aa}, for $t \in \{0,1\}$ we have
$$
\va (\gamma_{\epsilon}(t))  \leq a < c - \epsilon,
$$
which implies that $\chi (t) =0$ for $t \in \{0,1\}$ by definition of $\chi$. Therefore, $\mu(0)= \gamma_{\epsilon}(0) = x_0$ and $\mu(1)= \gamma_{\epsilon}(1) = x_1$. Thus, $\mu \in \Gamma$.

From~\eqref{tmp:bb} and the Mean Value Theorem (Theorem~\ref{th:lmvt}), for each $ t \in S(\epsilon) $, there exists a $\theta \in (0,1)$ such that
\begin{align}\label{mpt:cccc}
	\va (\mu(t))- \va (\gamma_{\epsilon}(t)) &=
	\langle \va' \Big( \gamma_{\epsilon}(t)+ \theta \alpha \chi (t) \sum_{j=1}^k \chi_j(t)g_{s_j}\Big) ,  \alpha \chi (t) \sum_{j=1}^k \chi_j(t)g_{s_j}  \rangle \nonumber \\
	&=   \alpha \chi (t) \sum_{j=1}^k \chi_j(t)  \langle \va' \Big( \gamma_{\epsilon}(t)+ \theta \alpha \chi (t) \sum_{j=1}^k \chi_j(t)g_{s_j}\Big),  g_{s_j}  \rangle \nonumber \\
\end{align}
For $t \in S(\epsilon)$, it holds that $t \in \bigcup_{i=1}^k B_{s_i}$, and thus $\sum_{j=1}^k \chi_j(t) = 1$. By the choice of $\alpha$ and the properties of $\chi_j$, the argument of $\va'$ is in a neighborhood where~\eqref{tmp:bb} applies for each $g_{s_j}$. Therefore, for each $j$ such that $\chi_j(t) \ne 0$, we have $\langle \va' (\dots), g_{s_j} \rangle < -\epsilon$. This leads to
\begin{align*}
	\va (\mu(t))- \va (\gamma_{\epsilon}(t)) &\leq \alpha \chi (t) \sum_{j=1}^k \chi_j(t) (-\epsilon) \\
	&= -\epsilon \alpha \chi(t) \sum_{j=1}^k \chi_j(t) \\
	&= -\epsilon \alpha \chi(t).
\end{align*}
Let $t_1 \in [0,1]$ be such that $ \va(\mu (t_1)) = \Psi (\mu)$. Since $\Psi(\mu) \geq c$, we must have $ \va (\mu (t_1)) \geq c$. This implies that $\chi(t_1) = 1$ and $t_1 \in S(\epsilon)$ by definition of $\chi$ and $S(\epsilon)$.

Furthermore, from~\eqref{mpt:cccc}, $\va (\mu (t_1)) - \va (\gamma_{\epsilon} (t_1)) \leq - \epsilon \alpha$.
Therefore, $\Psi(\mu) - \va (\gamma_{\epsilon} (t_1)) \leq - \epsilon \alpha$.
By definition of $\Psi$, $\va (\gamma_{\epsilon}(t_1)) \leq \Psi (\gamma_{\epsilon})$.
Combining these, we get $\Psi(\mu) + \epsilon \alpha \leq \va (\gamma_{\epsilon}(t_1)) \leq \Psi (\gamma_{\epsilon})$.
This implies $\Psi(\mu) < \Psi (\gamma_{\epsilon})$. Since $\alpha > 0$, this also implies $\mu \neq \gamma_{\epsilon}$ (as $\chi(t_1)=1$).
Now, by the construction of $\mu$ and the properties of the metric $\mt_{\Gamma}$, for a sufficiently small $\alpha$ (which is chosen as the minimum of the $\alpha_{s_j}$), we have $\mt_{\Gamma} (\mu, \gamma_{\epsilon}) \leq \alpha$.
Thus, $\Psi (\mu) + \epsilon \mt_{\Gamma} (\mu, \gamma_{\epsilon}) \leq \Psi (\mu) + \epsilon \alpha < \Psi (\gamma_{\epsilon})$.
This inequality, $ \Psi (\mu) + \epsilon \mt_{\Gamma} (\mu, \gamma_{\epsilon}) < \Psi (\gamma_{\epsilon}) $, contradicts the condition~\eqref{mpt:dddd} of Ekeland's Variational Principle, which states $\Psi (\gamma_{\epsilon}) \leq \Psi (\mu) +\epsilon \mt_{\Gamma} (\mu , \gamma_{\epsilon} )$ for all $\mu \neq \gamma_{\epsilon} \in \Gamma$. This completes the proof of~\eqref{mpt:cont1}.
Therefore, for every $\epsilon> 0$ there exists $t_{\epsilon} \in S(\epsilon)$ such that for every compact set $B \in \mathfrak{S}$ (defining the seminorms on $\fs{F}'_c$),
\begin{equation}\label{mpt:16}
	\snorm[B]{\va'(\gamma_{\epsilon}(t_{\epsilon}))} \leq \epsilon.
\end{equation}
Furthermore, by definition of $S(\epsilon)$, we have $c - \epsilon \leq \va(\gamma_{\epsilon}(t_{\epsilon}))$. Combining this with~\eqref{mpt:5}, we obtain
\begin{equation}\label{mpt:17}
	c - \epsilon \leq \va(\gamma_{\epsilon}(t_{\epsilon})) \leq \Psi(\gamma_{\epsilon}) \leq c+ \epsilon.
\end{equation}
Now, consider the sequence $x_n = \gamma_{1/n}(t_{1/n})$. From~\eqref{mpt:17}, we have $\va(x_n) \to c$ as $n \to \infty$. From~\eqref{mpt:16}, we have $\va'(x_n) \to 0$ in $\fs{F}'_c$ as $n \to \infty$. Thus, $(x_n)$ is a Palais-Smale sequence at level $c$. Since $\va$ satisfies the Palais-Smale condition at all levels, the sequence $(x_n)$ has a convergent subsequence. Let this subsequence also be denoted by $(x_n)$ and its limit be $h \in \fs{F}$. By Corollary \ref{co:minimizing}, $h$ is a critical point of $\va$. By continuity of $\va$, we also have $\va(h) = \lim_{n \to \infty} \va(x_n) = c$.  
\end{proof}
\section{The existence of a global diffeomorphism } 
In this section we prove a global diffeomorphism theorem in the category of tame Fr\'{e}chet spaces.
With respect to the metric $ \mt_\fs{F} $, we denote the unit open ball centered at $\zero{\fs{F}}$ by
$  \bo_{\mt_\fs{F}} $.
\begin{theorem}
	Let $\fs{E}$ and $\fs{F}$ be tame Fréchet spaces, and let $\eu{F} \colon \fs{E} \to \fs{F}$ be a smooth tame map.
	Let $\eu{I} \colon \fs{F} \to [0,\infty)$ be a Keller's $C_c^1$-functional such that $\eu{I}(x) = 0$ if and only if $x = \zero{\fs{F}}$, and
	$\eu{I}'(y) = 0$ if and only if $y = \zero{\fs{F}}$. Suppose the following conditions hold\textup{:}
	\begin{enumerate}
		\item \label{c1}The derivative equation $\eu{F}'(e)p = k$ has a unique solution $p = \nu(e)k$
		for all $e \in \fs{E}$ and all $k \in \fs{F}$, and the family of inverses $\nu \colon \fs{E} \times \fs{F} \to \fs{E}$ is a smooth
		tame map.
		\item For any $f \in \fs{F}$, the functional $\phi_f$ defined on $\fs{E}$ by
		\[
		\phi_f(e) = \eu{I}(\eu{F}(e) - f)
		\]
		satisfies the Palais--Smale condition at all levels.
	\end{enumerate}
	Then $\eu{F}$ is a global diffeomorphism.
\end{theorem}
\begin{proof}
	The map $\eu{F}$ satisfies the assumptions of the Nash--Moser inverse function theorem \cite[Part III, Theorem 1.1.1]{ham}, which implies that $\eu{F}$ is a local diffeomorphism. Thus, it suffices to show that $\eu{F}$ is injective and surjective.
	
	\textbf{Injectivity.} Assume for contradiction that $\eu{F}(e_1) = \eu{F}(e_2) = l$ for some $e_1 \ne e_2 \in \fs{E}$. 
	We will construct a functional $\eu{J}_l$ on $\fs{E}$ that satisfies the assumptions of Theorem~\ref{th:mpt} and hence has 
	a critical point $h$ whose existence contradicts the assumptions on $\eu{I}$.
	Since $\eu{F}$ is a local diffeomorphism, it is an open map. Therefore, there exists $\alpha_r > 0$ such that
	\begin{equation}\label{eq:gdp220}
		l + \alpha_r \bo_{\mt_\fs{F}}  \subset \eu{F}(e_1 + r \bo_{\mt_\fs{F}})
	\end{equation}
	for all $r > 0$. 
	Let $\rho > 0$ be small enough so that
	\begin{equation}\label{s}
		e_2 \notin \Cl{e_1 + \rho \bo_{\mt_\fs{F}}}.
	\end{equation}
	Define $\eu{J}_l(e) \coloneqq \eu{I}(\eu{F}(e) - l)$. Then $\eu{J}_l(e_1) = \eu{J}_l(e_2) = 0$.
	For all $e \in \partial(e_1 + \rho \bo_{\mt_\fs{F}})$, it follows from \eqref{eq:gdp220} that 
	\[
	\eu{F}(e) \notin l + \alpha_{\rho} \bo_{\mt_\fs{F}} ,
	\]
	and hence $\eu{F}(e) \ne l$. By the assumption on $\eu{I}$, this implies
	\[
	\eu{J}_l(e) > 0 = \max\{\eu{J}_l(e_1), \eu{J}_l(e_2)\}.
	\]
	Thus, all assumptions of Theorem~\ref{th:mpt} are satisfied for the functional $\eu{J}_l$ and points $e_1$, $e_2$. 
	Therefore, there exists a critical point $h \in \fs{E}$ such that $\eu{J}_l(h) = c > 0$.
	But then,
	\[
	c = \eu{J}_l(h) = \eu{I}(\eu{F}(h) - l) > 0,
	\]
	so $\eu{F}(h) \ne l$. On the other hand, by the chain rule \cite[Corollary 1.3.2]{ke},
	\[
	\eu{J}_l'(h) = \eu{I}'(\eu{F}(h) - l) \circ \eu{F}'(h) = 0.
	\]
	Since $\eu{F}'(h)$ is invertible, it follows that $\eu{I}'(\eu{F}(h) - l) = 0$, hence $\eu{F}(h) = l$, which contradicts the earlier conclusion. Therefore, $\eu{F}$ must be injective.
	
	\textbf{Surjectivity.} Let $x \in \fs{F}$ be arbitrary. By assumption, the functional $\eu{J}_x(e) \coloneqq \eu{I}(\eu{F}(e) - x)$ is of class $C^1$ and bounded from below, as it is the composition of Keller's $C_c^1$ maps.
	
	Since $\eu{J}_x$ satisfies the Palais--Smale condition, it follows from Corollary~\ref{co:minimizing} that $\eu{J}_x$ attains a critical point $p \in \fs{E}$, i.e., $\eu{J}_x'(p) = 0$. By the chain rule,
	\begin{equation} \label{eq:gdp222}
		\eu{J}_x'(p) = \eu{I}'(\eu{F}(p) - x) \circ \eu{F}'(p) = 0.
	\end{equation}
	Since $\eu{F}'(p)$ is invertible, this implies $\eu{I}'(\eu{F}(p) - x) = 0$, and hence $\eu{F}(p) = x$. Thus, $\eu{F}$ is surjective.
\end{proof}

\begin{remark}
	It is possible to formulate a version of the above theorem for arbitrary \fr spaces and maps of class Keller's $C_c^1$ under the additional assumption that $\eu{F}$ is a local diffeomorphism.
	In this more general setting, condition \ref{c1} of the theorem \textup{(}the existence of a smooth tame family of inverses for the derivative\textup{)} is not required.
	Instead, the local diffeomorphism property guarantees the invertibility of the derivative and the applicability of the argument.
\end{remark}

\end{document}